\documentclass[english,10pt]{article}
\usepackage{pst-node}
\usepackage{amsmath}
\usepackage{pifont}
\usepackage{mathrsfs}
\usepackage{amssymb}
\usepackage{stmaryrd}
\usepackage{tipa}
\usepackage{tikz}
\usepackage{amssymb}
\usepackage{latexsym}
\usepackage{amsthm}
\usepackage{fancyhdr}
\usepackage{amsfonts}
\usepackage{graphicx}
\usepackage{verbatim}
\usepackage{color}

\newtheorem{theorem}{Theorem}[section]

\newtheorem{definition}[theorem]{Definition}
\newtheorem{conjecture}[theorem]{Conjecture}

\newtheorem{lemma}[theorem]{Lemma}

\newtheorem{proposition}[theorem]{Proposition}

\newcommand{\Z}{\mbox{$\mathbb Z$}}
\newcommand{\M}{\mbox{$\mathcal M$}}

\newcommand{\setZTHR}{\mbox{$\langle\Z_3\rangle$}}

\begin{document}
\title{Nowhere-zero $3$-flow and $\Z_3$-connectedness in Graphs with Four Edge-disjoint Spanning Trees}

\author{Miaomiao Han\thanks{Department of Mathematics, West Virginia University,
Morgantown, WV 26506, USA.
Email: mahan@mix.wvu.edu}, Hong-Jian Lai\thanks{Department of Mathematics, West Virginia University,
Morgantown, WV 26506, USA.
Email: hjlai@math.wvu.edu}, Jiaao Li\thanks{Department of Mathematics, West Virginia University,
Morgantown, WV 26506, USA.
Email: joli@mix.wvu.edu}}

\maketitle

\begin{abstract}
Given a zero-sum function $\beta : V(G) \rightarrow \mathbb{Z}_3$ with $\sum_{v\in V(G)}\beta(v)=0$,
an orientation $D$ of $G$ with $d^+_D(v)-d^-_D(v)= \beta(v)$ in $\mathbb{Z}_3$ for every vertex $v\in V(G)$
is called a $\beta$-orientation.
A graph $G$ is $\mathbb{Z}_3$-connected if $G$ admits a $\beta$-orientation for every
zero-sum function $\beta$.
Jaeger et al. conjectured that  every $5$-edge-connected graph is $\Z_3$-connected.
A graph is $\setZTHR$-extendable at vertex $v$ if any pre-orientation at
$v$ can be extended to a $\beta$-orientation of $G$ for any zero-sum function $\beta$.
We observe that if
every $5$-edge-connected essentially $6$-edge-connected graph is $\setZTHR$-extendable at any degree five vertex,
then the above mentioned conjecture by Jaeger et al. holds as well.
Furthermore,  applying the  partial flow extension method of Thomassen and of Lov\'{a}sz et al.,  we prove that every graph with at least 4 edge-disjoint spanning trees is $\Z_3$-connected.
Consequently, every $5$-edge-connected essentially $23$-edge-connected graph is
$\setZTHR$-extendable at degree five  vertex.
\end{abstract}

\section{Introduction}

We consider finite graphs without loops, but with possible multiple
edges, and follow \cite{BoMu08} for undefined
terms and notation. As in \cite{BoMu08},
$\kappa'(G)$ denotes the edge-connectivity of a graph $G$; and $d_D^+(v)$, $d_D^-(v)$ denote
the out-degree and the in-degree of a vertex in a digraph $D$, respectively.
Throughout this paper, $\Z$ denotes the set of integers, and
$A$ denotes an (additive) abelian group with identity $0$.  For an $m\in \mathbb{Z}$,
let $\mathbb{Z}_m$ be the set of integers modulo $m$, as well as the (additive)
cyclic group on $m$ elements. For vertex subsets $U,  W\subseteq V(G)$,
let $[U, W]_G = \{uw \in E(G)|u \in U, w \in W\}$; and for each $v \in V(G)$,
define $E_G(v)=[v, V(G)-v]_G$. The subscript $G$ may be omitted if $G$
is understood from the context.  An edge cut $X=[S,V(G)-S]$ in a connected graph
$G$ is {\bf essential} if at least two components of $G-X$ are nontrivial.
A graph is {\bf essentially $k$-edge-connected} if it does not
have an essential edge cut with fewer than $k$ edges.

For an integer $m>1$, a graph $G$ admits a {\bf mod $m$-orientation}
if $G$ has an
orientation $D$ such that at every vertex $v\in V(G)$,
$d^+_D(v)-d^-_D(v)\equiv 0 \pmod m$. Let $\M_m$ be the family of all
graphs admitting a mod $m$-orientation.
Let $k\ge 2$ be an integer and $G$ be a
graph with an orientation $D=D(G)$. For any vertex
$v\in V(G)$, let $E^+_D(v)$ denote the set of all edges
directed away from $v$, and let $E^-_D(v)$ denote the set of
all edges directed into $v$. A function $f: E(G)\rightarrow \{\pm 1, \pm2,\dots, \pm (k-1)\}$ is called a {\bf nowhere-zero $k$-flow} if
\[
\sum\limits_{e\in E^+_D(v)}f(e) ~-\sum\limits_{e\in E^-_D(v)}f(e)=0,
\mbox{ for any vertex $v\in V(G)$.}
\]
The well-known $3$-Flow Conjecture of Tutte is stated below.

\begin{conjecture}\label{tutte3flow}
  (Tutte   \cite{Tutt54}) Every $4$-edge-connected graph has a nowhere-zero $3$-flow.
\end{conjecture}
Tutte \cite{Tutt66} (see also Brylawski \cite{Bryl72},
Arrowsmith and Jaeger \cite{ArJa82}) indicated
that a graph $G$ has a nowhere-zero $k$-flow if and only if  $G$ has a nowhere-zero $\Z_k$-flow.
Moreover,  a graph has a nowhere-zero $3$-flow if and only if $G$ has a mod $3$-orientation (i.e. $G\in \M_3$).

Jaeger et al. \cite{JLPT92} introduced the notion of $\mathbb{Z}_k$-connectedness as a generalization of nowhere-zero flows. In this paper, we mainly focus on $\mathbb{Z}_3$-connectedness of graphs. A function $\beta : V(G) \rightarrow \mathbb{Z}_3$ is a zero-sum function of $G$ if $\sum_{v\in V(G)}\beta(v)=0$ in $\Z_3$. Let $Z(G, \Z_3)$ be the set of all  zero-sum functions of $G$.
An orientation $D$ of $G$ with $d^+_D(v)-d^-_D(v)= \beta(v)$ in $\mathbb{Z}_3$ for every vertex $v\in V(G)$
is called a {\bf $\beta$-orientation}. A mod $3$-orientation of $G$ is a $\beta$-orientation with $\beta(v)=0$ for every vertex $v\in V(G)$.
A graph $G$ is {\bf $\mathbb{Z}_3$-connected}
if, for every $\beta\in {Z}(G, \mathbb{Z}_{3})$, there is an orientation
$D$ such that $d^+_D(G)-d^-_D(G)\equiv \beta(v)  \pmod {3}$ for every vertex
$v\in V(G)$.  The collection of all $\Z_3$-connected graphs is denoted by $\langle\Z_3\rangle$. Jaeger et al. \cite{JLPT92}
proposed the following  Conjecture.

\begin{conjecture}\label{Jaegerz3}
  (Jaeger, Linial, Payan and Tarsi \cite{JLPT92}) Every $5$-edge-connected graph is $\Z_3$-connected.
\end{conjecture}

A graph $G$ with $z_0 \in V(G)$ is {\bf $\M_3$-extendable at vertex $z_0$}
if,  for any pre-orientation $D_{z_0}$ of $E_G(z_0)$ with $d_{D_{z_0}}^+(z_0) \equiv d_{D_{z_0}}^-(z_0) \pmod 3$,
$D_{z_0}$ can be extended to a mod $3$-orientation $D$ of $G$.
Kochol \cite{Koch01} showed that Conjecture
\ref{Jaegerz3} implies  Conjecture \ref{tutte3flow}.

\begin{theorem}\label{kochol}(Kochol \cite{Koch01}) The following  are equivalent.
\\
(i) Every $4$-edge-connected graph has a nowhere-zero $3$-flow.
\\
(ii) Every $5$-edge-connected graph has a nowhere-zero $3$-flow.
\\
(iii) Every $5$-edge-connected essentially $6$-edge-connected  graph is $\M_3$-extendable at every degree $5$ vertex.
\\
(iv) Every $4$-edge-connected graph with each vertex of degree $4$ or $5$ is $\M_3$-extendable at every vertex.
\end{theorem}

A graph is called {\bf $\setZTHR$-extendable at vertex $z_0$},
if,  for any $\beta\in Z(G, \Z_3)$ and any pre-orientation $D_{z_0}$ of $E_G({z_0})$ with $d_{D_{z_0}}^+(z_0)-d_{D_{z_0}}^-(z_0)\equiv\beta (z_0)\pmod 3$,
$D_{z_0}$ can be extended to a $\beta$-orientation $D$ of $G$.
In the next section, we shall prove  the following proposition on extendability at vertex $z_0$.

\begin{proposition}\label{extendingiff} Let $G$ be a graph and $z_0 \in V(G)$ be a vertex.
\\
(i) $G$ is $\setZTHR$-extendable at vertex $z_0$ if and only if $G-z_0$ is $\Z_3$-connected.
\\
(ii) If $G$ is $\setZTHR$-extendable at $z_0$, then $G$ is $\Z_3$-connected.
\end{proposition}

Thomassen \cite{Thom12} and Lov\'{a}sz et al. \cite{LTWZ13} utilized partial
flow extensions to obtain breakthroughs in $\Z_3$-connectedness and
modulo orientation problems.
Lov\'{a}sz, Thomassen, Wu and Zhang \cite{LTWZ13,Wuyz12} proved that every
$6$-edge-connected graph is $\Z_3$-connected. In fact, they have proved a stronger result.

\begin{theorem}\label{LTWZWu} (Lov\'{a}sz, Thomassen, Wu and Zhang \cite{LTWZ13} and Wu \cite{Wuyz12})
 Every $6$-edge-connected  graph is $\setZTHR$-extendable at any vertex of degree at most $7$.
\end{theorem}

{ Analogous to Theorem \ref{kochol}(iii) of Kochol, it is natural to suggest the following strengthening of Conjecture \ref{Jaegerz3}, which eliminates nontrivial $5$-edge-cut,
and whose truth would imply Conjecture \ref{Jaegerz3}, as to be shown in Section 3 of this paper.}

\begin{conjecture}\label{z3extending}
Every $5$-edge-connected essentially $6$-edge-connected graph is $\setZTHR$-extendable at any vertex of degree $5$.
\end{conjecture}

The main results of this paper are the following.
\begin{theorem}\label{mainthm}
  Every graph with $4$ edge-disjoint spanning trees is $\Z_3$-connected.
\end{theorem}

In response to Theorem \ref{kochol}(iii) of Kochol and providing some supporting evidence to Conjecture \ref{z3extending},  we  obtain a
partial result as stated below.
\begin{theorem}\label{ess23con} Each of the following holds.\\
(a)Every $5$-edge-connected essentially $23$-edge-connected graph is $\M_3$-extendable at any degree five vertex.\\
(b)Every $5$-edge-connected essentially $23$-edge-connected graph is $\setZTHR$-extendable at any degree five vertex.
\end{theorem}

Theorems \ref{mainthm} and \ref{ess23con} are immediate corollaries of a technical theorem, stated below as Theorem \ref{4treez3}, which would be proved via utilizing a method of Thomassen \cite{Thom12} and Lov\'{a}sz et al. in \cite{LTWZ13}.

Following Catlin \cite{Catl92}, let  $F(G,k)$ denote
the minimum number of additional edges that must be  added to
$G$ to result in a supergraph $G'$ of $G$ that has $k$ edge-disjoint spanning trees.
In particular, $G$ has $k$ edge-disjoint spanning trees if and only if
$F(G,k) = 0$.
It is known (\cite{WLYZ14,LaLL17}) that if $G$ is $\Z_3$-connected, then it contains two edge-disjoint spanning trees (i.e.
$F(G,2) = 0$). A cut-edge is called a {\bf bridge}.
The following provides a sufficient condition for graphs to be $\Z_3$-connected through number of edge-disjoint spanning trees.

\begin{theorem}\label{4treez3} Let $G$ be a graph.
\\
(i) Suppose that $F(G,4)\le 3$. Then $G$ is $\Z_3$-connected, unless $G$ contains a bridge. (Thus, $G$ is $\Z_3$-connected if and only if $\kappa'(G)\ge 2$. )
\\
(ii) Suppose that  $F(G,4)=0$. Then for any vertex $v \in V(G)$ with $d_G(v) \le 7$, if $\kappa'(G-v) \ge 2$, then
$G$ is $\setZTHR$-extendable at  $v$.
\end{theorem}

Prerequisites will be presented in the next section. In Section 3, we will study the relationship
among
Conjectures  \ref{tutte3flow},  \ref{Jaegerz3} and \ref{z3extending}.
 Theorems \ref{4treez3}, \ref{mainthm} and \ref{ess23con} will be proved in a subsequent section.

\section{Prerequisites}

In this section, we will justify Proposition \ref{extendingiff} and present other preliminaries.
For a graph $G$ and a vertex $z \in V(G)$, define
$N_G(z)   =\{v\in V(G): zv\in E(G)\}$.
For notation convenience, the algebraic manipulations in the proof of Proposition \ref{extendingiff} will be
over $\Z_3$.
\\

\noindent
{\bf Proof of Proposition \ref{extendingiff} } As Part (ii) is straightforward, we only prove Part (i).
Suppose that a graph  $G$ is $\setZTHR$-extendable at vertex $z_0$.
Let $D_{z_0}$ be a fixed pre-orientation of $E_G({z_0})$. We also use $D_{z_0}$ to denote the digraph
induced by the oriented edges of $D_{z_0}$. Define
\begin{equation} \label{b}
\mbox{ $b(v)=d_{D_{z_0}}^+(v)-d_{D_{z_0}}^-(v)$ for each $v\in N_G(z_0)\cup \{z_0\}$.}\end{equation}
Then $b(z_0) + \sum_{v \in N_G(z_0)} b(v) = 0$.

We are to prove  $G-z_0$ is $\mathbb{Z}_3$-connected. For any $\beta\in \mathbb{Z} (G-{z_0}, \mathbb{Z}_3)$, define
$$\beta'(v)=
\left\{\begin{array}{lll}
\beta(v)+b(v), \mbox{if~} v\in N_G(z_0); \\
b(z_0), \mbox{if~} v=z_0;\\
\beta(v),  \mbox{~othewise.}
\end{array}
\right.
$$
Then $\sum_{v \in V(G)} \beta'(v) = \sum_{v \in V(G-z_0)} \beta(v) + (b(z_0) + \sum_{v \in N_G(z_0)} b(v)) = 0$,
and so
$\beta'\in  {Z} (G, \mathbb{Z}_3)$. Since $G$ is $\setZTHR$-extendable at vertex $z_0$,
there exists an orientation $D'$ of $G$ such that $d^+_{D'}(v)-d^-_{D'}(v)=\beta'(v)$ for any vertex $v\in V(G)$
and $D'$ agrees with $D_{z_0}$ on $E_G(z_0)$.
Let $D$ be the restriction of $D'$ on $G-z_0$.
By the definition of $\beta'$, we have $d^+_D(v)-d^-_D(v)=\beta(v)$
for any vertex $v\in V(G-z_0)$, and so $G-z_0$ is $\mathbb{Z}_3$-connected.

Conversely, assume that $G-z_0$ is $\mathbb{Z}_3$-connected.
Let $\beta'\in {Z} (G, \mathbb{Z}_3)$, and $D_{z_0}$ be a pre-orientation of
$E_G({z_0})$ with $d^+_{D_{z_0}}(z_0)-d^-_{D_{z_0}}(z_0)=\beta'(z_0)$.
Define $b(v)$ as in (\ref{b}), and
$$\beta(v)=
\left\{\begin{array}{ll}
\beta'(v)-b(v), \mbox{if~} v\in N_G(z_0); \\
\beta'(v), \mbox{otherwise.}
\end{array}
\right.
$$
As $\sum_{v \in V(G - z_0)} \beta(v) = \sum_{v \in V(G)} \beta'(v) = 0$,
we have $\beta\in  {Z} (G-z_0, \mathbb{Z}_3)$. Since $G-z_0 \in \setZTHR$,
there exists an orientation $D'$ of $G-z_0$ satisfying  $d^+_{D'}(v)-d^-_{D'}(v)=\beta'(v)$,
for any vertex $v\in V(G-z_0)$.
Combine $D'$ and $D_{z_0}$ to obtain an orientation $D$ of $G$. Then for any vertex $v\in V(G)$,
depending on $v = z_0$ or not, we always have $d^+_D(v)-d^-_D(v)=\beta'(v)$, and so $G$ is $\setZTHR$-extendable at vertex $z_0$.
This completes the proof of Proposition  \ref{extendingiff}.
\qed

\vskip 0.3cm

Let $G$ be a graph and $\beta\in Z(G,\Z_3)$. Define an integer valued mapping
$\tau : 2^{V(G)}  \mapsto \{0,\pm 1, \pm 2, \pm 3\}$  as follows:  for each vertex $x\in V(G)$,
\begin{eqnarray}\nonumber
  \tau (x) \equiv \left\{\begin{array}{ll}
\beta(x)\pmod 3; \\
d(x) \pmod 2.
\end{array}
\right.
\end{eqnarray}

For a vertex set $A\subset V(G)$, denote $\beta(A)\equiv \sum_{v\in A} \beta(v)\pmod 3$, $d(A)=|[A,V(G)-A]|$
and define
$\tau(A)$ to be
\begin{eqnarray}\nonumber
  \tau (A) \equiv \left\{\begin{array}{ll}
\beta(A)\pmod 3; \\
d(A) \pmod 2.
\end{array}
\right.
\end{eqnarray}

\begin{theorem}\label{partialextending}
(Lov\'{a}sz, Thomassen,  Wu and Zhang, Theorem 3.1 of \cite{LTWZ13})
Let $G$ be a graph, $\beta\in Z(G,\Z_3)$ and $z_0 \in V(G)$.
If $D_{z_0}$ is a pre-orientation of $E_G({z_0})$, and if
\\
(i) $|V(G)|\ge 3$,
\\
(ii) $d(z_0)\le 4 + |\tau(z_0)|$ and $d^+(z_0)-d^-(z_0)\equiv \beta(z_0) \pmod 3$, and
\\
(iii) $d(A)\ge 4+ |\tau(A)|$ for each nonempty $A\subseteq V(G)-\{z_0\}$ with $|V(G)- A|\ge 2$,
\\
then $D_{z_0}$ can be extended to a $\beta$-orientation of the entire graph $G$.
\end{theorem}

The following is an application of Theorem \ref{partialextending}.

\begin{lemma}\label{delete3edges}
Let $G$ be a $6$-edge-connected graph. Each of the following holds.
\\
(i) If $v \in V(G)$ with $d(v) \le 7$, then $G-v \in \setZTHR$.
\\
(ii) If $E_1 \subset E(G)$ with $|E_1|\le 3$, then $G-E_1 \in \setZTHR$.
\end{lemma}

\proof ~ (i) we may assume that $d_G(v)=7$ to prove the lemma. Otherwise, pick an edge $e \in E_G(v)$ and add an edge parallel to $e$, which results in still a $6$-edge-connected graph.
Take an arbitrary $\beta'\in {Z}(G-{v}, \mathbb{Z}_3)$. We shall show that $G-v$ has a $\beta'$-orientation.
Define $\beta(v)=3$.
We shall apply Theorem \ref{partialextending} by viewing $v$ as $z_0$ in Theorem \ref{partialextending}.
Since $d(v) = 7$, we have $|\tau(v)|=3$, and thus we can
orient  the edges $E_G(v)$ with an orientation $D_v$
so that $d_{D_v}^+(v) = 5$ and $d_{D_v}^-(v) = 2$.
Define $b(x)=d_{D_{v}}^+(x)-d_{D_{v}}^-(x)$ for each $x\in N_G(v)$ and set
\begin{equation} \label{beta-1}
\beta(x)=
\left\{\begin{array}{lll}
\beta'(x)+b(x), \mbox{if~} x\in N_G(v); \\
\beta(v), \mbox{if~} x=v;\\
\beta'(x),  \mbox{~othewise.}
\end{array}
\right.
\end{equation}
Then $\beta\in  {Z} (G, \mathbb{Z}_3)$. As $\kappa'(G) \ge 6$,
conditions (i)-(iii) of Theorem \ref{partialextending} are satisfied, and so by Theorem \ref{partialextending},
$G$ has a $\beta$-orientation $D$.
Let $D'$ be the restriction of $D$ on $G-v$. By (\ref{beta-1}),
$D'$ is a $\beta'$-orientation of $G-v$. This proves (i).

(ii) Since $\Z_3$-connectedness is preserved under adding edges, we may assume
that $|E_1|=3$. In graph $G$, subdivide each edge in $E_1$ with internal vertices $z_1,z_2,z_3$,
respectively. Identify $z_1,z_2,z_3$ to form a new vertex $z_0$ in the
resulted graph $G'$. By the construction of $G'$, we have $\kappa'(G') \ge 6$.
By Lemma \ref{delete3edges} (i), $G-E_1=G'-z_0 \in \setZTHR$.
\qed

\vspace{0.3cm}

For an edge set $X \subseteq E(G)$, the
{\bf contraction} $G/X$ is the graph obtained from $G$ by
identifying the two ends of each edge in $X$, and then deleting the
resulting loops. If $H$ is a subgraph of $G$, then we use $G/H$ for
$G/E(H)$.  For a vertex set $W \subset V(G)$ such that $G[W]$ is connected, we also use $G/W$ for $G/G[W]$.

\begin{lemma}\label{cf} (Proposition 2.1 of \cite{LaiH03})
Let $G$ be a graph.  Each of the following holds.
\\
(i) If $G\in \langle \mathbb{Z}_3\rangle$ and $e\in E(G)$, then $G/e \in \langle \mathbb{Z}_3\rangle.$
\\
(ii) If $H\subseteq G$ and if $H, G/H \in \langle \mathbb{Z}_3\rangle$, then $G\in \langle \mathbb{Z}_3\rangle.$
\end{lemma}

\section{Relationship among the conjectures}

A graph is called {\bf $\setZTHR$-reduced} if it does not have any nontrivial
$\Z_3$-connected subgraphs. By definition, $K_1$ is $\setZTHR$-reduced. The potential minimal counterexamples of Conjectures \ref{tutte3flow} and \ref{Jaegerz3} must be $\setZTHR$-reduced graphs.
As an example, it is routine to verify
that  the $4$-edge-connected non $\Z_3$-connected graph $J$ constructed
by Jaeger et al.\cite{JLPT92} (see Figure \ref{Jaegergraph}) is indeed a $\setZTHR$-reduced graph.
Applying Theorem \ref{partialextending}, we obtain the following.

\begin{lemma}\label{min5}
  Every $\setZTHR$-reduced graph has minimal degree at most $5$.
\end{lemma}
\proof
Suppose, to the contrary, that there is  a  $\setZTHR$-reduced graph $G$
with $\delta(G)\ge 6$. As a cycle of length $2$ is $\Z_3$-connected,
$G$ has no parallel edges and $|V(G)|\ge 4$. If $\kappa'(G)\ge 6$, then $G$ is
$\Z_3$-connected by Theorem \ref{LTWZWu}, contradicting to $G$ is  a
$\setZTHR$-reduced graph. For a vertex subset $W \subset V(G)$, let
$W^c = V(G) - W$. Among all those edge-cuts $[W, W^c]$
of size at most $5$ in $G$,  choose the  one with $|W|$ minimized. Let $v_c$ denote the vertex onto which $W^c$
is contracted in $G/W^c$. Obtain a graph $G'$ from $G/W^c$ by adding $6-d_{G/W^c}(v_c)$ edges
between $W$ and $v_c$. Then $\kappa'(G')\ge 6$ by the choice of $W$.
By Lemma \ref{delete3edges} (i), $G[W]= G'- v_c$ is $\Z_3$-connected, a contradiction. \qed
\\

We believe that the following strengthening of Lemma \ref{min5} holds as well, whose truth
implies Conjecture \ref{Jaegerz3}, as will be shown below in Proposition \ref{d5extending}.

\begin{conjecture}\label{atmost4}
Every $\setZTHR$-reduced graph has minimal degree at most $4$.
\end{conjecture}

\begin{proposition}\label{d5extending} Each of the following holds.
\\
(i)
Conjecture \ref{z3extending} implies Conjecture \ref{atmost4}.
\\
(ii)  Conjecture \ref{atmost4} implies Conjecture \ref{Jaegerz3}.
\end{proposition}

\proof { We shall prove (ii) first.} Assume that Conjecture \ref{atmost4} holds. Then
by the validity of Conjecture \ref{atmost4}, every graph with minimum degree at least 5
is not  $\setZTHR$-reduced.
Let $G$ be a counterexample to Conjecture \ref{Jaegerz3} with $|V(G)|$ minimized.
Since $\delta(G) \ge \kappa'(G) \ge 5$, $G$ is not  $\setZTHR$-reduced, and so
$G$ contains a nontrivial $\Z_3$-connected subgraph $H$. Since $\kappa'(G/H) \ge \kappa'(G) \ge 5$,
and since $|V(G)| > |V(G/H)|$, the minimality
of $G$ implies that $G/H$ is $\Z_3$-connected. By Lemma \ref{cf} (ii), $G$ must  be $\Z_3$-connected as well, contrary to
the assumption that $G$ is a counterexample of Conjecture \ref{Jaegerz3}.
This proves { (ii)}.

To prove (i), we use arguments similar to those in the proof of Lemma \ref{min5}.
By contradiction, we assume that Conjecture \ref{z3extending} holds
but there is a counterexample $G$ to Conjecture \ref{atmost4} with $|V(G)|$ minimized and with $\delta(G) \ge 5$.
By the validity of Conjecture \ref{z3extending},
$G$ must have an essential edge-cut of size at most $5$.
Among all those essential edge-cuts $[W, W^c]$
of size at most $5$,  choose the  one with $|W|$ minimized.
Let $v_c$ denote the vertex onto which $W^c$
is contracted in $G/W^c$. Adding some edges between $W$ and $v_c$ such that $v_c$ has degree $5$ in the new graph, and we still denote it $G/W^c$. Then we have $|W|\ge 2$,
and the minimality of $|W|$ forces that $G/W^c$ is an essentially $6$-edge-connected graph. By the assumption that Conjecture \ref{z3extending} holds,
$G/W^c$ is $\setZTHR$-extendable at $v_c$.
By Proposition \ref{extendingiff}, $G[W]=G/W^c - v_c \in \setZTHR$, contradicting to that $G$ is $\setZTHR$-reduced.
\qed
\\

In the rest of this section, we study the relationship between $\setZTHR$-extendability and edge deletions.
Theorem \ref{eqdeletetwoedges}
below indicates that deleting one or two adjacent edges does not make Conjecture
\ref{Jaegerz3} stronger. Theorem \ref{eqdeleteanytwoedges} and Proposition \ref{eqdeletethreeedges} below also describe the strength of Conjecture \ref{atmost4} and Conjecture \ref{z3extending} via edge deletions.

\begin{theorem} \label{eqdeletetwoedges} The following statements are equivalent.
\\
(i) Every $5$-edge-connected graph is $\Z_3$-connected.
\\
(ii) Every $5$-edge-connected graph deleting two  adjacent edges is $\Z_3$-connected.
\end{theorem}

\begin{theorem} \label{eqdeleteanytwoedges} The following statements are equivalent.
\\
(i) Every $\setZTHR$-reduced graph has minimal degree at most $4$.
\\
(ii) Every $5$-edge-connected graph deleting  any two  edges is $\Z_3$-connected.
\end{theorem}

\begin{proposition} \label{eqdeletethreeedges} The following statements are equivalent.
\\
(i) Every $5$-edge-connected essentially $6$-edge-connected graph is $\setZTHR$-extendable at any vertex of degree $5$.
\\
(ii) Every $5$-edge-connected  graph is $\setZTHR$-extendable at any vertex of degree $5$.
\\
(iii) Every $5$-edge-connected  graph deleting three {incident} edges of a degree $5$ vertex is $\Z_3$-connected.
\end{proposition}

 We shall justify Theorem
\ref{eqdeletetwoedges} and Theorem \ref{eqdeleteanytwoedges} by utilizing Kochol's method in \cite{Koch01}. In \cite{Koch01}, Kochol applies $\M_3$-extension on a degree $5$ vertex and converts it into degree $3$ vertices, which helps him establish Theorem \ref{kochol}. Unlike mod $3$-orientations, direct application of the method above does not seem to help on $\setZTHR$-extension for certain $\beta$-orientation. We observe that some edge deletions behave similarly as extension, as showed in Proposition \ref{extendingiff} and the theorems above. This is part of the reason why we would like to prove Theorem \ref{4treez3} in the form of edge deletions.

 A lemma is needed to
prove Theorems  \ref{eqdeletetwoedges} and \ref{eqdeleteanytwoedges}.

\begin{definition} \label{e-sum}
Let $G_1$ be a  graph with $e =u_1v_1 \in E(G_1)$,
and $G_2(u_2, v_2)$ be a graph with distinguished (and distinct) vertices of $u_2, v_2$.
Let $G_1 \oplus_e G_2$ be a graph obtained from the disjoint union
of $G_1 - e$ and $G_2$ by identifying $u_1$ and $u_2$ to form a vertex $u$, and by identifying
$v_1$ and $v_2$ to form a vertex $v$. Thus for $i \in \{1, 2\}$, we can view $u = u_i$ and $v = v_i$ in $G_i$.
Note that even if  $e$ and $u_2, v_2$ are given, $G_1 \oplus_e G_2$ may not be unique. Thus we use
$G_1 \oplus_e G_2$ to denote any one of the resulting graph.
\end{definition}

\begin{lemma}\label{2sum}
Let $G_1$ and $G_2$ be nontrivial graphs with $e \in E(G_1)$.

(i) If $G_1$ and $G_2$ are not $\Z_3$-connected graphs, then $G_1 \oplus_e G_2$ is not $\Z_3$-connected.

(ii) If $G_1$ and $G_2$ are $\setZTHR$-reduced graphs, then $G_1 \oplus_e G_2$ is a $\setZTHR$-reduced graph..
\end{lemma}

\proof ~(i) The proof is similar to those of Lemma 1 in \cite{Koch01} and of Lemma 2.5 in \cite{DeXY06}.
Let $G = G_1 \oplus_e G_2$. We shall adopt the notation in Definition \ref{e-sum}.
Fix $ i \in \{1, 2\}$.
Since $G_i$ is not $\Z_3$-connected,
there exists a $\beta_i \in Z(G_i, \Z_3)$ such that
$G_i$ does not have a $\beta_i$-orientation.
Define $\beta: V(G) \mapsto \Z_3$ as follows:

$$\beta(x)=
\left\{\begin{array}{lll}
\beta_1(x), \mbox{if~} x\in V(G_1)-\{u_1,v_1\}; \\
\beta_2(x), \mbox{if~} x\in V(G_2)-\{u_2,v_2\}; \\
\beta_1(x)+\beta_2(x),  \mbox{if~} x\in \{u,v\}.
\end{array}
\right.
$$
As $\sum_{z \in V(G)} \beta(z) =  \sum_{i=1}^2 \sum_{z \in V(G_i)} \beta_i(z)$,
we have $\beta\in Z(G, \Z_3)$. It remains to show $G$ does not have a $\beta$-orientation.
By contradiction, assume that $G$ has a $\beta$-orientation $D$.
Let $D_2$ be the restriction of $D$ on $E(G_2)$. Then $d^+_{D_2}(x)-d^-_{D_2}(x)=\beta_2(x)$
in $\Z_3$ for any $x\in V(G_2)-\{u_2,v_2\}$.
Since $G_2$ does not have a $\beta_2$-orientation,
we must have $d^+_{D_2}(u)-d^-_{D_2}(u)\neq \beta_2(u)$ in $\Z_3$. Thus, we have either
\begin{eqnarray} \label{+1}
  d^+_{D_2}(u)-d^-_{D_2}(u)= \beta_2(u)+1~~\text{and}~~ d^+_{D_2}(v)-d^-_{D_2}(v) = \beta_2(v)-1,
\end{eqnarray}
or
\begin{eqnarray} \label{-1}
  d^+_{D_2}(u)-d^-_{D_2}(u)= \beta_2(u)-1~~\text{and}~~ d^+_{D_2}(v)-d^-_{D_2}(v) = \beta_2(v)+1.
\end{eqnarray}

Let $D_1'$ be the restriction of $D$ on $E(G_1)-e$. If (\ref{+1}) holds,
then both $d^+_{D_1'}(u)-d^-_{D_1'}(u)= \beta_1(u)-1$ and $d^+_{D_1'}(v)-d^-_{D_1'}(v) = \beta_1(v)+1$.
Obtain an orientation $D_1$ of $G_1$ from $D_1'$ by orienting $e=u_1v_1$ from $u_1$ to $v_1$.
If (\ref{-1}) holds,
then both  $d^+_{D_1'}(u)-d^-_{D_1'}(u)= \beta_1(u)+1$ and $d^+_{D_1'}(v)-d^-_{D_1'}(v) = \beta_1(v)-1$.
Obtain an orientation $D_1$ of $G_1$ from $D_1'$ by orienting $e=u_1v_1$ from $v_1$ to $u_1$.
In either case, $D_1$ is a $\beta_1$-orientation of $G_1$, contrary to
the choice of $\beta_1$.  (ii) follows from (i) by the definition of $\setZTHR$-reduced graph. This proves the lemma.
\qed

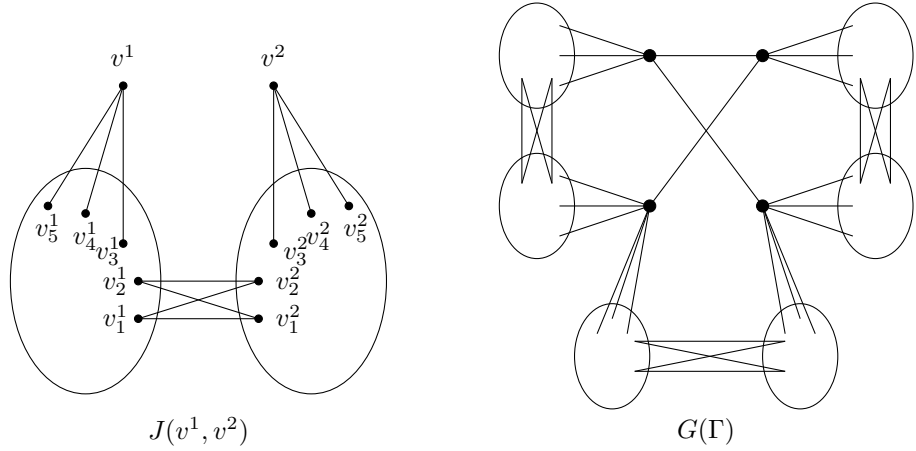
\begin{figure}
\begin{center}
\begin{tikzpicture}
\draw (0,0)
ellipse (1cm and 1.5cm);
\draw (3,0)
ellipse (1cm and 1.5cm);
\filldraw[black] (0.5,2.6)  circle (0.05cm);
\filldraw[black] (-0.5,1) circle (0.05cm);
\filldraw[black] (0,0.9) circle (0.05cm);
\filldraw[black] (0.5,0.5) circle (0.05cm);
\filldraw[black] (0.7,0) circle (0.05cm);
\filldraw[black] (0.7,-0.5) circle (0.05cm);
\filldraw[black] (2.5,2.6) circle (0.05cm);
\filldraw[black] (3.5,1) circle (0.05cm);
\filldraw[black] (3,0.9) circle (0.05cm);
\filldraw[black] (2.5,0.5) circle (0.05cm);
\filldraw[black] (2.3,0) circle (0.05cm);
\filldraw[black] (2.3,-0.5) circle (0.05cm);
\draw [-] (0.7,0)--(2.3,0);
\draw [-] (0.7,0)--(2.3,-0.5);
\draw [-] (0.7,-0.5)--(2.3,0);
\draw [-] (0.7,-0.5)--(2.3,-0.5);
\draw [-] (0.5,2.6)--(-0.5,1);
\draw [-] (0.5,2.6)--(0,0.9);
\draw [-] (0.5,2.6)--(0.5,0.5);
\draw [-] (2.5,2.6)--(3.5,1);
\draw [-] (2.5,2.6)--(3,0.9);
\draw [-] (2.5,2.6)--(2.5,0.5);
\node at (0.5,3){$v^1$};
\node at (2.5,3){$v^2$};
\node at (-0.5,0.7){$v^1_5$};
\node at (0,0.6){$v^1_4$};
\node at (0.3,0.4){$v^1_3$};
\node at (0.4,0){$v^1_2$};
\node at (0.4,-0.5){$v^1_1$};
\node at (3.6,0.7){$v^2_5$};
\node at (3.1,0.6){$v^2_4$};
\node at (2.8,0.4){$v^2_3$};
\node at (2.7,0){$v^2_2$};
\node at (2.7,-0.5){$v^2_1$};

\draw (6,3)
ellipse (0.5cm and 0.7cm);
\draw (10.5,3)
ellipse (0.5cm and 0.7cm);
\draw (6,1)
ellipse (0.5cm and 0.7cm);
\draw (7,-1)
ellipse (0.5cm and 0.7cm);
\draw (10.5,1)
ellipse (0.5cm and 0.7cm);
\draw (9.5,-1)
ellipse (0.5cm and 0.7cm);
\filldraw[black] (7.5,3) circle (0.08cm);
\filldraw[black] (9,3) circle (0.08cm);
\filldraw[black] (7.5,1) circle (0.08cm);
\filldraw[black] (9,1) circle (0.08cm);
\draw [-] (7.5,3)--(9,3);
\draw [-] (7.5,3)--(9,1);
\draw [-] (7.5,1)--(9,3);
\draw [-] (7.5,3)--(6.3,3.4);
\draw [-] (7.5,3)--(6.3,3);
\draw [-] (7.5,3)--(6.3,2.6);
\draw [-] (7.5,1)--(6.3,1.4);
\draw [-] (7.5,1)--(6.3,1);
\draw [-] (7.5,1)--(6.3,.6);
\draw [-] (7.5,1)--(7.2,-0.7);
\draw [-] (7.5,1)--(7,-0.5);
\draw [-] (7.5,1)--(6.8,-0.7);
\draw [-] (9,3)--(10.2,3.4);
\draw [-] (9,3)--(10.2,3);
\draw [-] (9,3)--(10.2,2.6);
\draw [-] (9,1)--(10.2,1.4);
\draw [-] (9,1)--(10.2,1);
\draw [-] (9,1)--(10.2,0.6);
\draw [-] (9,1)--(9.7,-0.7);
\draw [-] (9,1)--(9.5,-0.5);
\draw [-] (9,1)--(9.3,-0.7);
\draw [-] (5.8,2.7)--(5.8,1.3);
\draw [-] (5.8,2.7)--(6.2,1.3);
\draw [-] (6.2,2.7)--(6.2,1.3);
\draw [-] (6.2,2.7)--(5.8,1.3);
\draw [-] (10.3,2.7)--(10.3,1.3);
\draw [-] (10.3,2.7)--(10.7,1.3);
\draw [-] (10.7,2.7)--(10.7,1.3);
\draw [-] (10.7,2.7)--(10.3,1.3);
\draw [-] (7.3,-0.8)--(9.3,-1.2);
\draw [-] (7.3,-1.2)--(9.3,-0.8);
\draw [-] (9.3,-0.8)--(7.3,-0.8);
\draw [-] (9.3,-1.2)--(7.3,-1.2);
\node at (1.5,-2.0){$J(v^1, v^2)$};
\node at (8.25,-2.0){$G(\Gamma)$};
\end{tikzpicture}
\end{center}
\caption{ The  construction in Theorem \ref{eqdeletetwoedges}}\label{construction}
\end{figure}

~

\noindent
{\bf Proof of Theorem \ref{eqdeletetwoedges}.} It suffices to prove that (i) implies (ii). By contradiction, assume that (i) holds
and that there exists a graph $\Gamma$ with $\kappa'(\Gamma) \ge 5$ and with two distinct
adjacent edges $vv_1, vv_2 \in E(\Gamma)$, where $v_1$ and $v_2$ may or maynot be
distnict, such that $\Gamma - \{vv_1, vv_2\} \notin \setZTHR$. As $\kappa'(\Gamma) \ge 5$,
$|E_{\Gamma}(v)| \ge 5$.
Let $K \cong K_4$ with $V(K) = \{w_1, w_2, w_3, w_4\}$.

We assume first that $v_1 \neq v_2$ in $\Gamma$, and
use $L(v_1, v_2)$ to denote  $\Gamma - \{vv_1, vv_2\}$ with $v_1$ and
$v_2$ being two distinguished vertices.
For $1 \le j \le 2$, let $\phi_j: L_j(v_1^j, v_2^j) \mapsto L(v_1, v_2)$ be
a graph isomorphism with $\phi_j(v^j) = v$, $\phi_j(v_1^j) = v_1$  and  $\phi_j(v_2^j) = v_2$.
Define $J(v^1, v^2) = K \oplus_{w_1w_2} L_1(v_1^1, v_2^1)  \oplus_{w_3w_4} L_2(v_1^2, v_2^2)$.
Let  $J^k(v^1, v^2)$, ($1 \le k \le 3$), be three isomorphic copies of  $J(v^1, v^2)$,
and define $G(\Gamma) = K  \oplus_{w_1w_2} J^1(v^1, v^2) \oplus_{w_2w_3} J^2(v^1, v^2) \oplus_{w_3w_4} J^3(v^1, v^2)$,
as depicted in  Figure \ref{construction}. By the definition of $G(\Gamma)$, $G(\Gamma)$
contains six subgraphs $H_i$, ($1 \le i \le 6$), each of which is isomorphic to $\Gamma - \{vv_1, vv_2\}$.

It is known that $K \notin \setZTHR$. As $\Gamma - \{vv_1, vv_2\} \notin \setZTHR$, it follows
from Lemma \ref{2sum} that  $J(v^1, v^2) \notin \setZTHR$, and so by repeated applications of
Lemma \ref{2sum}, $G(\Gamma)  \notin \setZTHR$.

Let $W \subseteq E(\Gamma)$ be a minimum edge cut of $G(\Gamma)$.  If for any $i$, $|W \cap E(H_i)| = 0$, then
$W$ is an edge cut of the graph $G(\Gamma)/(\cup_{i=1}^6 H_i)$, and so it is straightforward to check that $|W| \ge 5$.
Hence we assume that for some $i$, $W \cap E(H_i) \neq \emptyset$. Then
$\Gamma - \{vv_1, vv_2\}$ contains an edge subset $W_i'$ corresponding to $W \cap E(H_i)$
under the isomorphism between $\Gamma - \{vv_1, vv_2\}$ and $H_i$.
If $W_i'$ does not separate the neighbors of $v$ and $\{v_1, v_2\}$ in $\Gamma$,
then $W_i'$ is an edge cut of $\Gamma$, and so
$|W| \ge |W'| \ge \kappa'(\Gamma) \ge 5$. Hence by symmetry, we assume that
$v$ and $v_1$ are in different components of $\Gamma - W_i'$.
Since $\kappa'(\Gamma) \ge 5$, we have $|W_i'| \ge \kappa'(\Gamma - \{vv_1, vv_2\}) = 5-2 = 3$.
By the definition of $G(\Gamma)$, $G(\Gamma) - E(H_i)$ contains 2 edge-disjoint $(v, v_1)$-paths,
which implies that $|W - E(H_i)| \ge 2$, and so
$|W| = |W \cap E(H_i)| + |W - E(H_i)| \ge 3+2=5$.
We conclude that $\kappa'(G(\Gamma)) \ge 5$. By Theorem \ref{eqdeletetwoedges}(i),
we have $G(\Gamma) \in \setZTHR$, which leads to a contradiction to
the fact that $G(\Gamma)  \notin \setZTHR$.

Next we assume that $v_1 = v_2$.
Then for $j = 1, 2$, $v_1^j = v_2^j$ in $L_j(v_1^j, v_2^j)$. In this case,
we  differently define $J(v^1, v^2)$ to be the graph obtained from the disjoint union
of $L_1(v_1^1, v_2^1)$ and $L_2(v_1^2, v_2^2)$ by identifying
$v_1^1$ with $v^2_1$. Since $L_1(v_1^1, v_2^1)$ is a block of $J(v^1, v^2)$,
$J(v^1, v^2) \notin \setZTHR$. We again define
$G(\Gamma) = K  \oplus_{w_1w_2} J^1(v^1, v^2) \oplus_{w_2w_3} J^2(v^1, v^2) \oplus_{w_3w_4} J^3(v^1, v^2)$.
Then by  Lemma \ref{2sum}, $G(\Gamma)  \notin \setZTHR$. By a similar argument as shown above,
we again conclude that  $\kappa'(G(\Gamma)) \ge 5$, and so by Theorem \ref{eqdeletetwoedges}(i),
$G(\Gamma) \in \setZTHR$. This contradiction establishes the theorem.
\qed

\begin{figure}
\begin{center}
\begin{tikzpicture}

\filldraw[black] (0,0) circle (0.05cm);
\filldraw[black] (0,1) circle (0.05cm);
\filldraw[black] (1,0) circle (0.05cm);
\filldraw[black] (1,1) circle (0.05cm);
\draw [-] (0,0)--(0,1);\draw [-] (0,0)--(1,0);\draw [-] (0,0)--(1,1);
\draw [-] (0,1)--(1,0);\draw [-] (0,1)--(1,1);\draw [-] (1,0)--(1,1);
\node at (-.3, .1){$x_1$};\node at (-.3,1.1){$x_2$};\node at (1.3,.1){$x_3$};\node at (1.3,1.1){$x_4$};
\node at (3.7, .1){$x_{11}$};\node at (3.7,1.1){$x_{12}$};\node at (5.3,.1){$x_{9}$};\node at (5.3,1.1){$x_{10}$};
\node at (2, 2.6){$x_7$};\node at (2,4.2){$x_5$};\node at (3,2.6){$x_8$};\node at (3,4.2){$x_6$};

\filldraw[black] (4,0) circle (0.05cm);
\filldraw[black] (4,1) circle (0.05cm);
\filldraw[black] (5,0) circle (0.05cm);
\filldraw[black] (5,1) circle (0.05cm);
\draw [-] (4,0)--(4,1);\draw [-] (4,0)--(5,0);\draw [-] (4,0)--(5,1);
\draw [-] (4,1)--(5,0);\draw [-] (4,1)--(5,1);\draw [-] (5,0)--(5,1);

\filldraw[black] (2,2.93) circle (0.05cm);
\filldraw[black] (2,3.93) circle (0.05cm);
\filldraw[black] (3,2.93) circle (0.05cm);
\filldraw[black] (3,3.93) circle (0.05cm);
\draw [-] (2,2.93)--(3,2.93);\draw [-] (2,2.93)--(3,3.93);\draw [-] (2,2.93)--(2,3.93);
\draw [-] (2,3.93)--(3,2.93);\draw [-] (2,3.93)--(3,3.93);\draw [-] (3,2.93)--(3,3.93);


\draw plot [smooth,tension=1.5] coordinates{(0,1)(1,3)(2,3.93)};
\draw plot [smooth,tension=1.5] coordinates{(1,1)(1.5,2.2)(2,2.93)};
\draw plot [smooth,tension=1.5] coordinates{(5,1)(4,3)(3,3.93)};
\draw plot [smooth,tension=1.5] coordinates{(4,1)(3.5,2.2)(3,2.93)};
\draw plot [smooth,tension=1.5] coordinates{(1,0)(2.5,-0.3)(4,0)};
\draw plot [smooth,tension=1.5] coordinates{(0,0)(2.5,-1)(5,0)};
\end{tikzpicture}
\end{center}
\caption{the graph $J$ : a $4$-edge-connected $\setZTHR$-reduced graph}\label{figureJaeger} \label{Jaegergraph}
\end{figure}
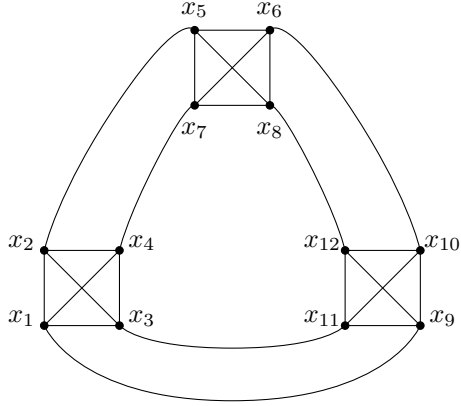

~

We need the following splitting theorem of Mader\cite{Made78} before proceeding the next proof. For two distinct vertices $x, y$, let $\lambda_G(x, y)$ be the maximum number of edge-disjoint paths connecting $x$ and $y$ in $G$. The following Mader's theorem
asserts that  local edge-connectivity is preserved under splitting.

\begin{theorem}\label{maderthm}(Mader \cite{Made78})
  Let $G$ be a graph and let $z$ be a non-separating vertex of $G$ with degree at least $4$ and $|N_G(z)|\ge 2$. Then there exist two edges $v_1z, v_2z$ in $G$ such that, splitting $v_1z, v_2z$, the resulting graph $G'=G-v_1z - v_2z + v_1v_2$ satisfies $\lambda_{G'}(x, y)=\lambda_G(x, y)$ for any two vertices $x, y$ different from $z$.
\end{theorem}

\begin{figure}
  \begin{center}
\begin{tikzpicture}
\node at (0, -3) {$H(w_3^1, w_3^2)$};
\node at (7.5, -3) {$G^*$};
\draw (0,0)
ellipse (1cm and 1.5cm);
\filldraw[black] (0.45,1.1) circle (0.05cm);\node at (.1, 1.1) {$u_1$};
\filldraw[black] (0.65,0.5) circle (0.05cm);\node at (.3, 0.5) {$u_2$};
\filldraw[black] (2,1.2) circle (0.05cm); \node at (2, 1.5) {$w_3^1$};
\draw [-] (0.45,1.1)--(2,1.2);
\draw [-] (0.65,0.5)--(2,1.2);

\filldraw[black] (0.65,-0.2) circle (0.05cm);\node at (.3, -0.2) {$v_1$};
\filldraw[black] (0.55,-0.8) circle (0.05cm);\node at (.2, -0.8) {$v_2$};
\filldraw[black] (2,-0.65) circle (0.05cm);\node at (2, -.95) {$w_3^2$};
\draw [-] (0.65,-0.2)--(2,-0.65);
\draw [-] (0.55,-0.8)--(2,-0.65);

\filldraw[black] (4,0) circle (0.05cm);
\filldraw[black] (4,-1) circle (0.05cm);
\filldraw[black] (5,0) circle (0.05cm);
\filldraw[black] (5,-1) circle (0.05cm);
\draw [-] (4,0)--(5,-1);\draw [-] (4,0)--(5,0);
\draw [-] (4,-1)--(5,0);\draw [-] (4,-1)--(5,-1);
\draw (3.25,-0.5)
ellipse (0.35cm and 0.75cm);
\draw [-] (4,0)--(3.3,0.1);
\draw [-] (4,0)--(3.3,-0.3);
\draw [-] (4,-1)--(3.3,-0.7);
\draw [-] (4,-1)--(3.3,-1.1);
\draw (5.75,-0.5)
ellipse (0.35cm and 0.75cm);
\draw [-] (5,0)--(5.7,0.1);
\draw [-] (5,0)--(5.7,-0.3);
\draw [-] (5,-1)--(5.7,-0.7);
\draw [-] (5,-1)--(5.7,-1.1);

\filldraw[black] (10,0) circle (0.05cm);
\filldraw[black] (10,-1) circle (0.05cm);
\filldraw[black] (11,0) circle (0.05cm);
\filldraw[black] (11,-1) circle (0.05cm);
\draw [-] (10,0)--(11,-1);\draw [-] (10,0)--(11,0);
\draw [-] (10,-1)--(11,0);\draw [-] (10,-1)--(11,-1);
\draw(9.25,-0.5)
ellipse (0.35cm and 0.75cm);
\draw [-] (10,0)--(9.3,0.1);
\draw [-] (10,0)--(9.3,-0.3);
\draw [-] (10,-1)--(9.3,-0.7);
\draw [-] (10,-1)--(9.3,-1.1);
\draw(11.75,-0.5)
ellipse (0.35cm and 0.75cm);
\draw [-] (11,0)--(11.7,0.1);
\draw [-] (11,0)--(11.7,-0.3);
\draw [-] (11,-1)--(11.7,-0.7);
\draw [-] (11,-1)--(11.7,-1.1);
\draw [-] (11,0)--(11.7,0.1);
\draw [-] (11,0)--(11.7,-0.3);
\draw [-] (11,-1)--(11.7,-0.7);
\draw [-] (11,-1)--(11.7,-1.1);


\filldraw[black] (7,3) circle (0.05cm);
\filldraw[black] (7,2) circle (0.05cm);
\filldraw[black] (8,3) circle (0.05cm);
\filldraw[black] (8,2) circle (0.05cm);
\draw [-] (7,3)--(8,2);\draw [-] (7,3)--(7,2);
\draw [-] (7,2)--(8,3);\draw [-] (8,3)--(8,2);
\draw(7.5,1.25)
ellipse (0.75cm and 0.35cm);
\draw [-] (7,2)--(7.15,1.3);
\draw [-] (7,2)--(6.8,1.3);
\draw [-] (8,2)--(7.8,1.3);
\draw [-] (8,2)--(8.15,1.3);
\draw(7.5,3.75)
ellipse (0.75cm and 0.35cm);
\draw [-] (7,3)--(7.15,3.7);
\draw [-] (7,3)--(6.8,3.7);
\draw [-] (8,3)--(7.8,3.7);
\draw [-] (8,3)--(8.15,3.7);

\draw plot [smooth,tension=1.5] coordinates{(4,0)(5.5,2)(7,3)};
\draw plot [smooth,tension=1.5] coordinates{(5,0)(6,1.5)(7,2)};
\draw plot [smooth,tension=1.5] coordinates{(8,2)(9,1.5)(10,0)};
\draw plot [smooth,tension=1.5] coordinates{(11,0)(9.5,2)(8,3)};
\draw plot [smooth,tension=1.5] coordinates{(4,-1)(7.5,-2.5)(11,-1)};
\draw plot [smooth,tension=1.5] coordinates{(5,-1)(7.5,-1.9)(10,-1)};

\end{tikzpicture}
\end{center}
\caption{The construction in Theorem \ref{eqdeleteanytwoedges}}
\end{figure}
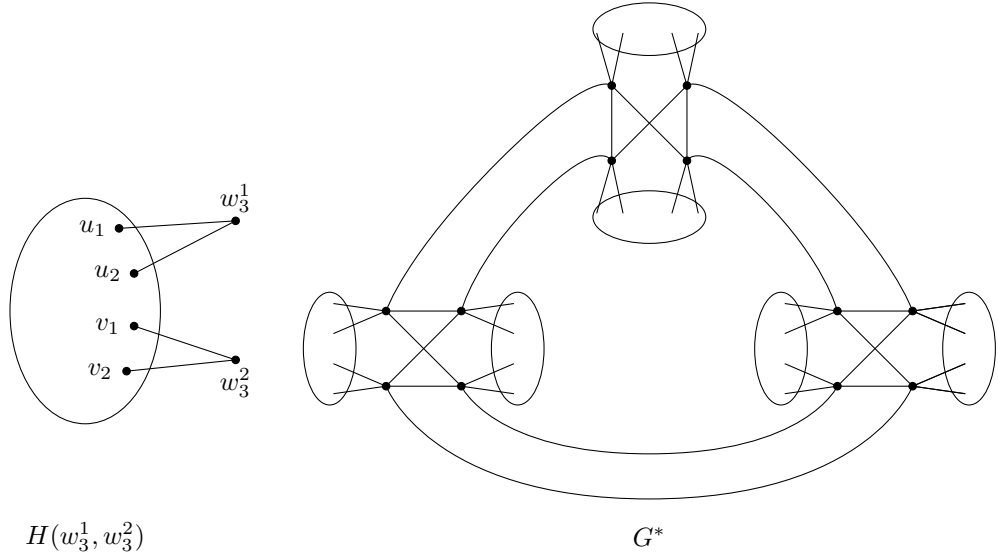


~

\noindent
{\bf Proof of Theorem \ref{eqdeleteanytwoedges}. } (i) $\Rightarrow$ (ii). By contradiction, assume that (i) holds
and that there exists a $5$-edge-connected graph $\Gamma$ with $|V(\Gamma)|$ minimized and with two distinct
 edges $u_1u_2, v_1v_2 \in E(\Gamma)$, where $u_1$ and $v_1$ may or maynot be
distinct, such that $G =\Gamma - \{u_1u_2, v_1v_2\} \notin \setZTHR$.
By the minimality of  $|V(\Gamma)|$,
$G$ must be a $\setZTHR$-reduced graph.
For $i=1, 2$, let $K^i \cong K_3$ with $V(K^i) = \{w_1^i, w_2^i, w_3^i\}$.
Define $K(v_1, v_2) = K^1 \oplus_{w_1^1w_2^1} G(u_1, u_2)$ and
$H(w_3^1, w_3^2)= K^2 \oplus_{w_1^2w_2^2}K(v_1, v_2)$. As $K_3$ and $G$
are $\setZTHR$-reduced graphs,  by Lemma \ref{2sum}(ii), $H(w_3^1, w_3^2)$ is a $\setZTHR$-reduced graph.
Moreover, $H(w_3^1, w_3^2)$ has exactly two vertices of degree 2, namely $w_3^1, w_3^2$, and the other
vertices of $H(w_3^1, w_3^2)$ have degree at least $5$.

Let $J$ be the graph as depicted in Figure \ref{Jaegergraph}
with $V(J) = \{x_1, \dots, x_{12}\}$. Obtain a graph $G^*$ by attaching copies of $H(w_3^1, w_3^2)$
and applying $\oplus_{e}$ operation for each $e=x_{2i-1}x_{2i}$, $1\le i\le 6$,
as depicted in Figure 3. Then we have $\delta(G^*) \ge 5$.
By the validity of (i), $G^*$ is not $\setZTHR$-reduced.
On the other hand, as $K_3$ and $G$ are $\setZTHR$-reduced, it follows by
Lemma \ref{2sum}(ii) that $H(w_3^1, w_3^2)$ is also $\setZTHR$-reduced. As $J$ and
$H(w_3^1, w_3^2)$ are $\setZTHR$-reduced, we conclude by
Lemma \ref{2sum}(ii) that $G^*$ is also $\setZTHR$-reduced, contrary to the fact
that $G^*$ is not $\setZTHR$-reduced, as implied by (i). This shows that
(i) implies (ii).

(ii) $\Rightarrow$ (i). Assume that (ii) holds. Then (ii) implies that every $5$-edge-connected graph is $\Z_3$-connected.
Let $G$ be a counterexample to (i). Then $G$ is a $\setZTHR$-reduced graph with $\delta(G)\ge 5$.
If $\kappa'(G) \ge 5$, then by (ii), $G$ itself is $\Z_3$-connected, contrary
to the assumption that  $G$ is $\setZTHR$-reduced.
Hence $\kappa'(G) \le 4$. Since $\delta(G) \ge 5$,
$G$ must have an essential edge-cut of size at most $4$.
Among all essential edge-cuts $[W, W^c]$
of size at most $4$,  choose one with $|W|$ minimized. Since $G$ is a $\setZTHR$-reduced graph,
$G[W]$ is also a  $\setZTHR$-reduced graph. Moreover, it is possible to add
two new  edges to $G[W]$ to result in a $5$-edge-connected graph.
If $|[W, W^c]|\le 3$, we obtain a graph $G[W]^+$ from $G[W]$
by appropriately adding two new edges (possibly parallel) joining
vertices in $W$ so that $\delta(G[W]^+) \ge 5$, and so
by the minimality of $|W|$, we have  $\kappa'(G[W]^+) \ge 5$.
By the validity of (ii), we conclude that $G[W]$
is $\Z_3$-connected. Since $\delta(G) \ge 5$, $G[W]$ is a nontrivial subgraph of $G$.
This contradicts the assumption that $G$ is a $\setZTHR$-reduced graph.

Hence we assume that $|[W, W^c]|=4$. Let $H=G/{W^c}$ and $z$ be the vertex onto which $G[W^c]$
is contracted, and denote $E_H(z)= \{e_1, e_2, e_3, e_4\}$ with $e_i = zv_i$, $1 \le i \le 4$.
Since $E_H(z)$ may contain parallel edges, the $v_i$'s do not have to be distinct.
By the minimality of $W$ and Menger's theorem,
we have $\lambda_H(x, y)\ge 5$ for  any two vertices $x, y\in V(H)-\{z\}$.

Suppose first that $H[E_H(z)]$
contains parallel edges. Assume that $z$ and $v_1$ are joined by at least 2 edges.
Define $H'' = H/H[\{z, v_1\}]$. By the minimality of $W$, we have $\kappa'(H'') \ge 5$.
As $|E_H(z) - E(H[\{z, v_1\}])| \le 2$, it follows
by (ii) that $G[W] = H'' - (E_H(z) - E(H[\{z, v_1\}]))$
is $\Z_3$-connected, contrary to the assumption that
$G$ is $\setZTHR$-reduced.

Hence we assume that $H[E_H(z)]$
contains no parallel edges, and so the $v_i$'s are 4 distinct vertices.
By Theorem \ref{maderthm}, we may assume that
the graph $H'=H -v_1z -v_2z+ v_1v_2$ satisfies $\lambda_{H'}(x, y)=\lambda_H(x, y)\ge 5$
for any two vertices $x, y\in V(H')-\{z\}$. This implies that the graph $H'' = H'/\{zv_3\}$
is $5$-edge-connected. By (ii), $G[W] \cong H'' - \{v_1v_2, e_4\} \in \setZTHR$,
contrary to the assumption that
$G$ is $\setZTHR$-reduced.
\qed
\\

Proposition \ref{eqdeletethreeedges} indicates certain implications of
Conjecture \ref{z3extending}.
The proof of Proposition \ref{eqdeletethreeedges} is similar to that of Proposition \ref{d5extending}
and is omitted.

\section{Proofs of Theorems \ref{mainthm}, \ref{ess23con} and \ref{4treez3}}

Theorems \ref{mainthm}, \ref{ess23con} and \ref{4treez3} will be proved in this section. We start with a lemma.

\begin{lemma}\label{lifting}
  Let $G$ be a graph, $v$ be a vertex of $G$ with degree at least $4$  and $vv_1, vv_2\in E_G(v)$. If $G_1=G-v+v_1v_2$ is $\Z_3$-connected, then $G$ is $\Z_3$-connected.
\end{lemma}
\proof~  Let $G_2=G-vv_1- vv_2+v_1v_2$. As $|[v,V(G)-v]_{G_2}| = d_G(v) - 2 \ge 2$, we have $G_2/G_1\in \setZTHR$.
Since $G_1\in\setZTHR$ and $G_2/G_1\in \setZTHR$, it follows by Lemma \ref{cf} that $G_2\in \setZTHR$.
By Lemma 2.4 of \cite{LaiH03},
 $G_2  \in \setZTHR$ implies that $G\in \setZTHR$.
\qed

For an integer $k > 0$, it is known (see \cite{Nash64}, or more explicitly,  Lemma 3.1 of \cite{LLLX10}
or Lemma 3.4 of \cite{LiLC09}) that if $F(H,k) > 0$ for any nontrivial proper subgraph $H$ of $G$,  then
\begin{equation} \label{FGk}
F(G,k) = k(|V(G)| - 1) - |E(G)|.
\end{equation}

\vspace{0.1cm}
\noindent
{\bf Proof of Theorem \ref{4treez3}.} Assume that Theorem \ref{4treez3} (i) holds and that $G$ is a graph with
$F(G,4) = 0$. If $v \in V(G)$ with $d_G(v) \le 7$ satisfies $\kappa'(G - v) \ge 2$, then $F(G-v, 4) \le 3$ and so by
Theorem \ref{4treez3} (i), $G - v$ is $\Z_3$-connected. It follows from Proposition \ref{extendingiff} that
$G$ is $\setZTHR$-extendable at vertex $v$. Thus if (i) holds, then (ii)
would follow as well. Hence
it suffices to show that
\begin{equation} \label{statement}
\mbox{ if $F(G,4)\le 3$ and $\kappa'(G) \ge 2$,  then
$G \in \setZTHR$.  }
\end{equation}
We argue by contradiction and assume that
\begin{equation} \label{ex}
\mbox{ $G$ is a counterexample to (\ref{statement}) with $|V(G)|+|E(G)|$ minimized.}
\end{equation}
As (i) holds if $|V(G)|\le 2$, we assume that $|V(G)| \ge 3$.  By assumption, there exists a set
$E_1$ of edges not in $G$ with $|E_1| = F(G,4)$
such that $G^+ = G + E_1$ contains four edge-disjoint spanning trees, denoted  $T_1, T_2, T_3, T_4$.
\\

\noindent{\bf Claim 1:} Each of the following holds.
\\
(i) For any nontrivial proper subgraph $H$ of $G$, $H \notin \setZTHR$  and
$F(H,4) \ge 3$.
\\
(ii) $G$ is $4$-edge-connected.

Let $H$ be a nontrivial proper subgraph of $G$.
As  $F(G/H, 4) \le 3$ (see, for example, Lemma 2.1 of  \cite{LiLC09}),
if $H \in \setZTHR$, then by (\ref{ex}) and $\kappa'(G/H)\ge\kappa'(G) \ge 2$, we have $G/H   \in \setZTHR$, and so by
Lemma \ref{cf},  $G \in \setZTHR$, contrary to (\ref{ex}). Hence
we must have $H \notin \setZTHR$.
If $F(H,4) \le 2$, then by $\kappa'(H)\ge 2$ and (\ref{ex}), we have
$H \in \setZTHR$, contrary to the fact that $H \notin \setZTHR$.
This proves Claim 1(i).

To prove Claim 1(ii), assume that $G$ has a minimum edge-cut $W$ with $|W| \le 3$.
Let $H_1$, $H_2$ be the two components of $G-W$. By (i) and by (\ref{FGk}),
we have
\[
F(H_1, 4) + F(H_2, 4) =  \sum_{i=1}^2 [4(|V(H_i)| - 1) - E(H_i)|] = F(G,4) - 4 + |W| \le |W|-1\le 2.
\]
This, together with the fact that $W$ is a minimum edge-cut,
implies that $\kappa'(H_i) \ge 2$  for
each $i \in \{1, 2\}$. Since $|V(G)|\ge 3$, at least one of $H_1$ and $H_2$ is nontrivial,
contrary  to Claim 1(i). Thus Claim 1(ii) must hold.
\\

\noindent{\bf Claim 2:} $E(G^+)=\cup_{i=1}^4E(T_i)$.

Suppose that there exists $e\in E(G^+)-\cup_{i=1}^4E(T_i)$. The minimality of $E_1$ indicates that $E_1 \subseteq \cup_{i=1}^4E(T_i)$, and thus $e\in E(G)$. Let $G' = G-e$. Then $G'$
is a spanning subgraph of $G$ with $F(G', 4) = F(G,4) \le 3$ and $\kappa'(G')\ge 3$ by Claim 1(ii). As $G' \in \setZTHR$
implies $G \in \setZTHR$, Claim 2 follows from (\ref{ex}).
\\

\noindent{\bf Claim 3:}  Each of the following holds.
\\
(i)  $G^+$ has no subgraph $H^+$ with $1 < |V(H^+)| < |V(G^+)|$ such that $F(H^+,4)=0$.
\\
(ii) $\kappa'(G^+) \ge 5$ and $G^+$ does not have an essentially $5$-edge-cut.
\\
(iii) $G^+$ has no vertex of degree $5$.

Argue by contradiction to show Claim 3(i) and choose
a subgraph $H^+$ of $G^+$ with $1 < |V(H^+)| < |V(G^+)|$ and $F(H^+,4)=0$
such that $|V(H^+)|$ minimized. By Claim 2, if $X =V(H^+)$, then
$H^+=G^+[X]$.  If $|X|=2$, then
by Claim 1(i), Claim 2 and $F(H^+,4)=0$, we conclude
that $E(G[X])$ consists of a cut edge of $G$, contrary to Claim 1(ii).
Hence we assume that $|X|\ge 3$.
Let $H = H^+ - E_1$. Then $H=G[X]$.
Since $F(H^+,4)=0$ and by Claim 2, $F(H,4) \le |E_1|= F(G,4) \le 3$. If $H$ has a cut edge $e$, then
by (\ref{FGk}) and as $|V(H)|\ge 3$, one component of $H-e$ must be nontrivial
and has 4 edge-disjoint spanning trees,
contrary to the minimality of $|V(H^+)|$. Hence $\kappa'(H) \ge 2$, and so by (\ref{ex}),
$H \in  \setZTHR$, contrary to Claim 1(i). This proves Claim 3(i).

If $W$ is a minimal 4-edge-cut or
an essential  $5$-edge-cut of $G^+$ with $G^+_1$ and $G^+_2$ being the two components of $G^+ - W$,
then by  (\ref{FGk}), there exists a nontrivial $H^+ \in \{G^+_1, G^+_2\}$
with $F(H^+,4)=0$,  contrary to Claim 3(i).  This proves Claim 3(ii).

We argue by contradiction to show Claim 3(iii).
Let $v_0$ be a vertex with  $d_{G^+}(v_0) = 5$,
$E_{G^+}(v_0)=\{e_1, e_2, e_3, e_4, e_5\}$, and $v_i$, $1 \le i \le 5$, be vertices
with $e_i = v_0v_i$. As $E_{G^+}(v_0)$ may contain parallel edges, the $v_i$'s are not necessarily distinct.
Since $F(G^+, 4)= 0$, we may assume that for $1 \le i \le 4$, $e_i \in E(T_i)$, and
$e_5 \in E(T_1)$. By Claim 1(ii), $|E_1 \cap  E_{G^+}(v_0)| \le 1$, and so we may assume
that $e_1 \in E(G)$.    By symmetry among $e_2, e_3, e_4$ and Claim 1(i) and (ii),
$e_1$ has at most one parallel edge, and thus we may assume $e_2\in E(G)$ and $v_2\neq v_1$.
Let $e_5''$ be an edge linking $v_1$ and $v_5$
but not in $E(G)$.
Define $G'' = G-v_0+{ v_1v_2}$ and
\[
E_1'' = \left\{
\begin{array}{ll}
E_1 & \mbox{ if $E_1 \cap E_{G^+}(v_0) = \emptyset$;}
\\
E_1 - E_{G^+}(v_0) & \mbox{ if $|E_1 \cap E_{G^+}(v_0)| =1$ and $e_5 \notin E_1$;}
\\
(E_1 - E_{G^+}(v_0)) \cup \{ e_5''\} & \mbox{ if $E_1 \cap E_{G^+}(v_0) = \{e_5 \}$.}
\end{array} \right.
\]
As for $i \in \{2,3, 4\}$, $T_i - v_0$ is a spanning tree of $G'' + E_1''$, and
$(T_1 - v_0) + e_5''$ is a spanning tree of $G''+E_1''$. It follows by
$|E_1''| \le |E_1| = 3$ that $F(G'', 4) \le 3$,  and $|V(G'')|+|E(G'')|<|V(G)|+|E(G)|$.
If $G''$ has a cut edge, then as $d_{G}(v_0) \le d_{G^+}(v_0) = 5$,
$G$ has a edge-cut $W'$ with $|W'| \le 3$, contrary to Claim 1(ii).
Thus $\kappa'(G'') \ge 2$.  By (\ref{ex}), $G''\in\setZTHR$.
Hence $G\in \setZTHR$ by Lemma \ref{lifting}, contrary to (\ref{ex}).
This proves Claim 3.

By Claim 3, $\kappa'(G^{+}) \ge 6$, and so
by Lemma \ref{delete3edges}(ii) and $F(G,4)\le 3$,
we have $G=G^+-E_1 \in \setZTHR$, contrary to (\ref{ex}). The proof is completed.
\qed

Theorem \ref{mainthm} is an immediate corollary of Theorem \ref{4treez3}, and we will prove Theorem \ref{ess23con} by a simple discharge argument.

The next lemma follows from arguments of Nash-Williams in \cite{Nash64}.
A detailed proof can be found in Theorem 2.4 of \cite{YaLL10}.

\begin{lemma} \label{gam}
Let $G$ be a nontrivial graph and let $k > 0$ be an integer.
If $|E(G)|\ge k(|V(G)|-1)$, then $G$ has a nontrivial
subgraph $H$ with $F(H,k) = 0$.
\end{lemma}

\noindent{\bf Proof of Theorem \ref{ess23con}.} It suffices to show (b). We shall show that every  $5$-edge-connected essentially $23$-edge-connected graph contains $4$ edge-disjoint spanning trees. Then Theorem \ref{ess23con}(b) follows from Theorem \ref{4treez3}(ii).

 Let $G$ be a counterexample  with $|E(G)|$ minimized.
Then $F(G,4)>0$ and  $|V(G)|\ge 4$.
If $|E(G)|\ge 4(|V(G)|-1)$, by Lemma \ref{gam}, there exists a non-trivial
subgraph $H$ with $F(H,4) = 0$.
By definition of contraction,  $G/H$ is $5$-edge-connected  and essentially $23$-edge-connected.
By the minimality of $G$, $G/H$ has $4$ edge-disjoint spanning trees.
As $H$ has $4$ edge-disjoint spanning trees, it follows that (see Lemma 2.1
of \cite{LiLC09}) $F(G,4) = 0$, contrary to the choice of $G$.
Hence we have
\begin{eqnarray}\label{4treesedges}
  |E(G)|<4(|V(G)|-1).
\end{eqnarray}
Since $|V(G)| \ge 4$ and  $G$ is essentially $23$-edge-connected,
for any edge $uv \in E(G)$, we have
\begin{eqnarray}\label{du}
d(u)+d(v)\ge 23+2.
\end{eqnarray}

For integers $i, k \ge 1$, define $D_i(G) = \{v \in V(G): d_G(v) = i\}$,  $D_{\le k}(G) = \cup_{i \le k} D_i(G)$,
and $D_{\ge k}(G) = \cup_{i \ge k} D_i(G)$. It follows from (\ref{du}) that $D_{\le 8}$ is an independent set.

Each vertex begins with charge equal to its degree. If $d(v)\ge 9$ and $vu\in E(G)$, then $v$ gives charge $\frac{d(v)-8}{d(v)}$ to $u$. Note that $G$ may contain parallel edges and the charge running through each edge adjacent to $v$.  Clearly, if $v\in D_{\ge 8}$, then $v$ will be left with charge $d(v)(1-\frac{d(v)-8}{d(v)})=8$.

For any vertex $x\in D_{\le 7}$, denote $d(x)=i\in \{5, 6, 7\}$. By (\ref{du}), $x$ will end with charge at least
$$i+\sum_{vx\in E(G)}\frac{d(v)-8}{d(v)}\ge i + \frac{25-i-8}{25-i}i=\frac{(42-2i)i}{25-i}\ge\min\{8, \frac{180}{19}, \frac{98}{9}\}=8,$$
a contradiction to (\ref{4treesedges}).
\qed

We remark that there exist $5$-edge-connected  and essentially $22$-edge-connected graphs do not contain $4$ edge-disjoint spanning trees. Lowing the constant $23$ may require new ideas and more elaborate work. As shown in Propositions \ref{d5extending} and \ref{eqdeletethreeedges}, lowing into $6$ would imply Conjectures \ref{tutte3flow} and \ref{Jaegerz3}.

\section{Two Applications}
Recall that a $\setZTHR$-reduced graph is a graph without nontrivial
$\Z_3$-connected subgraphs. The number of edges in a $\setZTHR$-reduced graph is often useful in reduction method and some inductive arguments. Theorem \ref{4treez3}, together with Lemma \ref{gam}, establishes an upper bound for the density of a $\setZTHR$-reduced graph.
\begin{lemma}
  Every $\setZTHR$-reduced graph on $n\ge 3$ vertices has at most $4n-8$ edges.
\end{lemma}

As defined in \cite{LLLM14}, a graph $G$ is {\bf strongly $\mathbb{Z}_{2s+1}$-connected} if, for
every $b : V(G) \rightarrow \mathbb{Z}_{2s+1}$ with $\sum_{v\in V(G)}b(v)=0$, there is an orientation $D$ such
that for every vertex $v\in V(G)$, $d^+_D(G)-d^+_D(G)\equiv b(v)  \pmod {2s+1}$.
Strongly $\mathbb{Z}_{2s+1}$-connected graphs are known as contractible configurations for
modulo $(2s+1)$-orientations. The following has recently been obtained.

\begin{proposition}\label{2strees}(\cite{LaLL17})
Every strongly $\mathbb{Z}_{2s+1}$-connected graph contains $2s$ edge-disjoint spanning trees.
\end{proposition}

By the monotonicity of circular flow (see, for example, \cite{GoTZ1998} or \cite{Zhan97}), it follows that every graph with a mod $5$-orientation also
has a mod $3$-orientation. It is not known, in general, whether a
strongly $\mathbb{Z}_{2k+3}$-connected graph is also strongly $\mathbb{Z}_{2k+1}$-connected.
As an application of  Proposition \ref{2strees}, if a graph $G$
is strongly $\mathbb{Z}_{5}$-connected graph, then $F(G,4) = 0$; it then follows from
Theorem \ref{4treez3} that $G \in \setZTHR$. Hence
we obtain the following proposition.

\begin{proposition}\label{sz3sz5}
   Every strongly $\mathbb{Z}_{5}$-connected graph is $\mathbb{Z}_{3}$-connected.
\end{proposition}

\end{document}